\documentclass[fleqn]{mat01}
\usepackage{times,mathtimy,amssymb}
\begin{document}

\setcounter{page}{205} \firstpage{205}

\newtheorem{theo}[defin]{\bf Theorem}
\newtheorem{coro}[defin]{\rm COROLLARY}
\newtheorem{lem}[defin]{Lemma}

\markboth{Rana Barua and Ashok Maitra}{Infinite products of Polish
spaces}

\title{Borel hierarchies in infinite products of Polish spaces}

\author{RANA BARUA and ASHOK MAITRA$^{*}$}

\address{Stat-Math Division, Indian Statistical Institute, Kolkata~700~108,
India\\
\noindent $^{*}$School of Statistics, University of Minnesota,
Minneapolis,
MN, USA\\
\noindent E-mail: rana@isical.ac.in; maitra@stat.umn.edu}

\volume{117}

\mon{May}

\parts{2}

\pubyear{2007}

\Date{MS received 22 August 2005; revised 28 August 2006}

\begin{abstract}
Let $H$ be a product of countably infinite number of copies of an
uncountable Polish space $X$. Let $\Sigma_\xi$ $(\overline
{\Sigma}_\xi)$ be the class of Borel sets of additive class $\xi$
for the product of copies of the discrete topology on $X$ (the
Polish topology on $X$), and let ${\cal B} = \cup_{\xi < \omega_1}
\overline{\Sigma}_\xi$. We prove in the L\'{e}vy--Solovay model
that
\begin{equation*}
\overline{\Sigma}_\xi =\Sigma_{\xi}\cap {\cal B}
\end{equation*}
for $1 \leq \xi < \omega_1$.
\end{abstract}


\maketitle

\section{Introduction}

Suppose $X$ is a Polish space and $N$ the set of positive
integers. We consider $H = X^N$ with two product topologies: (i)
the product of copies of the Polish topology on $X$, so that $H$
is again a Polish space and (ii) the product of copies of the
discrete topology on $X$. Define now the Borel hierarchy in the
larger topology on $H$. To do so, we need some notation. An
element of $H$ will be denoted by $h = (x_1, x_2,\dots,x_n,\dots)$
and for $m \in N$, $p_m(h)$ will denote the first $m$ coordinates,
that is, $p_m(h) = (x_1,x_2,\dots,x_m)$. For $n \in N$ and $A
\subseteq X^n$, $\hbox{cyl}(A)$ will denote the cylinder set with
base $A$, that is,
\begin{equation*}
\hbox{cyl}(A) = \{ h \in H\hbox{\rm :}\ p_n (h) \in A \}.
\end{equation*}
The Borel hierarchy for the larger topology on $H$ can now be
defined as follows:
\begin{equation*}
\Sigma_0 = \Pi_0 = \{ \hbox{cyl} (A)\hbox{\rm :}\ A \subseteq X^n,
\quad n \geq 1 \} \end{equation*} and for $\xi> 0$,
\begin{equation*}
\Sigma_\xi = \left(\bigcup_{\eta < \xi} \Pi_\eta\right)_\sigma,
\quad \Pi_\xi = \neg \Sigma_\xi.
\end{equation*}
The Borel hierarchy on $H$ with respect to the smaller topology is
defined in the usual way:
\begin{equation*}
\overline{\Sigma}_1 = \{ V\hbox{\rm :}\ V \ \mbox{is open in $H$
in the smaller topology}\},\ \overline{\Pi}_1 = \neg
\overline{\Sigma}_1
\end{equation*}
and, for $\xi > 1$,
\begin{equation*}
\overline{\Sigma}_\xi = \left(\bigcup_{\eta < \xi}
\overline{\Pi}_\eta\right)_\sigma; \;\; \overline{\Pi}_\xi = \neg
\overline{\Sigma}_\xi.
\end{equation*}
Let
\begin{equation*}
{\cal B} = \bigcup_{\xi < \omega_1} \overline{\Sigma}_\xi =
\bigcup_{\xi < \omega_1} \overline{\Pi}_\xi.
\end{equation*}
The problem we will address  in this article is whether
\begin{equation*}
\overline{\Sigma}_\xi =  \Sigma_\xi \cap {\cal B} \quad \mbox{for}
\ 1 \leq \xi < \omega_1.\tag{*}
\end{equation*}
To tackle the problem we will use the methods of effective
descriptive set theory. We therefore have to formulate the
lightface version of $(*)$. We refer the reader to [Mo] and [L1]
for definitions of lightface concepts. We take $X$ to be the
recursively presentable Polish space $\omega^{\omega}$ hereafter.

Define
\begin{equation*}
\Sigma^*_0 = \Pi^*_0 = \{\hbox{cyl} (A)\hbox{\rm :}\ A \ \mbox{is}
\ \Delta^1_1 \ \mbox{in} \ (\omega^{\omega})^n, \ n \geq 1 \},
\end{equation*}
and, for $\; 1 \leq \xi < \omega_1^{ck}$,
\begin{equation*}
\Sigma^*_\xi = \cup^1_1 (\cup_{\eta < \xi} \Pi^*_\eta)
\end{equation*}
and
\begin{equation*}
\Pi^*_\xi = \neg \Sigma^*_\xi,
\end{equation*}
where $\cup^1_1(\cup_{\eta < \xi} \Pi^*_\eta)$ is a $\Delta^1_1$
union of members of $\cup_{\eta < \xi} \Pi^*_\eta$. The lightface
analogue of $(*)$ is then
\begin{equation*}
\Sigma^*_\xi = \Delta^1_1 \cap \Sigma_\xi,\quad \mbox{for}\ 1 \leq
\xi < \omega_1^{ck}.\tag{**}
\end{equation*}
In order to state the  main result of the article, we equip
$\omega^{\omega}$ with the Gandy--Harrington topology, that is,
the topology whose base is the pointclass  of $\Sigma^1_1$ sets.
The key property of this topology is that it satisfies the Baire
category theorem (see [L1]). Consider now the following statement
of set theory:

\noindent ({\bf O})\quad Every subset of $\omega^{\omega}$ has the
Baire property with respect to the Gandy--Harrington topology.

The main result of the article can now be stated.

\begin{theo}[\!]
Assume {\rm ({\bf O})}. Let $1 \leq \xi < \omega_1^{ck}$. If $A$
and $B$ are $\Sigma^1_1$ subsets of $H$ such that $A$ can be
separated from $B$ by a $\Sigma_\xi$ set{\rm ,} then $A$ can be
separated from $B$ by a $\Sigma^*_\xi$ set.
\end{theo}

An immediate consequence is

\begin{coro}$\left.\right.$\vspace{.5pc}

\noindent {\rm ({\bf O})} implies $(^{**})$.
\end{coro}

The above results will be established in $\hbox{ZF}+\hbox{DC}$.
Maitra {\em et al} [Ma] proved ($^{*}$) for $\xi=1$ in
$\hbox{ZF}+\hbox{DC}$ by a boldface argument. We will provide a
lightface argument in the Appendix for ($^{**}$) when $\xi=1$.
Again this will be done in $\hbox{ZF}+\hbox{DC}$. Barua [Ba]
proved  Theorem~1.1 and Corollary~1.2. His proof was by induction
on $\xi$. However, he left out the proof of the base step
($\xi=1$). We will fill in the gap in this article. The proof of
Theorem 1.1 presented here parallels very closely that of Louveau
[L1], whereas the proof in [Ba] relies on the more abstract
developments of [L2]. In consequence, the proof given here is
somewhat simpler.

The paper is organized as follows. Section~2 is devoted to
definitions and notation. Section 3 contains the detailed proof of
Theorem~1.1 when $\xi=1$, while \S 4 sketches how the proof of
Theorem~1.1 can be completed by an inductive argument. In the
concluding section, we will prove ($^{*}$) under appropriate
hypotheses and also mention open problems.

\section{Definitions, notation and preliminaries}
\setcounter{defin}{0}

For $n \geq 1$, the Gandy--Harrington topology on
$(\omega^{\omega})^n$ will be denoted by $T^n$ and
the~Gandy--Harrington topology on $H$ will be denoted by
$T^\infty$. Following Louveau [L1], we define for each $\xi$ such
that $1 \leq \xi < \omega_1^{ck}$ a topology $T_\xi$ on $H$ having
for its base the pointclass $\Sigma^1_1 \cap \cup_{\eta <
\xi}\Pi_\eta$.

Let ${\cal S}$ be a second countable topology on
$(\omega^{\omega})^n$ (respectively, $H$). Let $A$ be a subset of
$(\omega^{\omega})^n$  (respectively, $H$). By the {\it
cosurrogate} of $A$ we mean the largest ${\cal S}$-open set $B$
such that $A \cap B$ is $T^n$-meager (respectively,
$T^\infty$-meager). The {\it surrogate} of $A$ is defined to be
the complement of the cosurrogate of $A$. When ${\cal S}$ is the
topology $T^n$, we denote  the surrogate (respectively,
cosurrogate) of $A$ by $\hbox{sur}^n(A)(\hbox{respectively,
cosur}^n(A))$. If $A \subseteq H$ and  ${\cal S}$ is the topology
$T_\xi$, the surrogate (respectively, cosurrogate) of $A$ will be
denoted by $\hbox{sur}_\xi(A)$ ($\hbox{respectively,
cosur}_\xi(A)$).

\begin{lem}
Let $m \geq 1$. If $A \subseteq (\omega^{\omega})^m$ is
$T^m$-open{\rm ,} then ${\rm sur}^m(A)$ is the $T^m$-closure of
$A$. Consequently{\rm ,} ${\rm sur}^m(A) - A$ is $T^m$-nowhere
dense.
\end{lem}

\begin{proof}
If $B$ is $\Sigma^1_1$ and $A \cap B$ is $T^m$-meager, then $A
\cap B$ must be empty, because $A \cap B$ is $T^m$-open and the
Baire category theorem holds for $T^m$. Consequently,
$\hbox{cosur}^m(A)$ is the union of basic open sets of the
$T^m$-topology which are disjoint with $A$. It follows that
$\hbox{sur}^m(A)$ is the $T^m$-closure of $A$. \hfill $\Box$
\end{proof}

\begin{lem}
Assume {\rm ({\bf O})}. Let $m \geq 1$. If $A \subseteq
(\omega^{\omega})^m${\rm ,} then $A \Delta \hbox{\rm sur}^m(A)$ is
$T^m$-meager.
\end{lem}

\begin{proof}
Observe that $\omega^{\omega}$ and $(\omega^{\omega})^m$ are
recursively isomorphic, so $(\omega^{\omega}, T^1)$ and
$((\omega^{\omega})^m, T^m)$ are homeomorphic. Hence it follows
from ({\bf O}) that  there is a $T^m$-open set $B$ such that $A
\Delta B$ is $T^m$-meager. So, if $D$ is a $\Sigma^1_1$ subset of
$(\omega^{\omega})^m$, then $A \cap D$ is $T^m$-meager iff $B \cap
D$ is $T^m$-meager, so that $\hbox{sur}^m(A) = \hbox{sur}^m(B)$.
Since $B$ is $T^m$-open, it follows from Lemma~2.1 that
$\hbox{sur}^m B-B$ is $T^m$-nowhere dense, hence $B \Delta
\hbox{sur}^m(B)$ is $T^m$- meager. Consequently, $A \Delta
\hbox{sur}^m(A)$ is $T^m$-meager. \hfill $\Box$
\end{proof}

Note that the converse of Lemma~2.2 is true. Indeed, if $A\Delta
\hbox{sur}^1(A)$ is $T^1$-meager for every $A\subseteq
\omega^{\omega}$, then, as is easy to verify, $A$ has the Baire
property with respect to $T^1$ for every $A\subseteq
\omega^{\omega}$, that is, ({\bf O}) holds.

\section{The case $\pmb{\xi = 1}$}

\setcounter{defin}{0}

In this section we will prove Theorem~1.1 when $\xi =
1$.\pagebreak

Following [L1], we fix a coding pair $(W,C)$ for the $\Delta^1_1$
subsets of $H$, that is,

\begin{enumerate}
\renewcommand\labelenumi{(\roman{enumi})}
\leftskip .3pc
\item $W$ is a $\Pi^1_1$ subset of $\omega$;
\item $C$ is a $\Pi^1_1$ subset of $\omega \times H$;
\item the relations `$n \in W$ \& $C(n,h)$' and `$n \in W$ \&
$\neg C(n,h)$' are both $\Pi^1_1$;
\item for every $\Delta^1_1$ subset $A$ of $H$, there is $n \in W$ such
that $A = C_n \stackrel{\rm def.}{=} \{ h \in H\hbox{\rm :}\
C(n,h)\}.$
\end{enumerate}

Define $W_0$ as follows:
\begin{align*}
&m \in W_0 \leftrightarrow m \in W\ \& \ (\exists n \geq
1)(\forall h)(\forall h')(C(n,h) \;\&\; p_n(h)\\[.5pc]
&\quad\, = p_n(h') \rightarrow C(n,h')).
\end{align*}
Then $W_0 $ is $\Pi^1_1$. Indeed, $W_0$ is just the set of codes
of $\Delta^1_1$ cylinder subsets of $H$.

\begin{lem}
If $A$ is a $\Sigma^1_1$ subset of $H${\rm ,} then ${\rm cl}_1(A)$
is $\Pi_1$ and $\Sigma^1_1${\rm ,} hence $T_2$-open{\rm ,} where
${\rm cl}_1(A)$ is the $T_1$-closure of $A$.
\end{lem}

\begin{proof}
Indeed, for any $A$, $\hbox{cl}_1(A)$ is $\Pi_1$, because it is a
countable intersection of $\Pi_1$ sets. Now suppose $A$ is
$\Sigma^1_1$. Then
\begin{align*}
h \notin \hbox{cl}_1(A) &\leftrightarrow (\exists n \geq
1)(\exists B)(B \ \mbox{is a } \ \Sigma^1_1 \; \mbox{subset of} \;
(\omega^{\omega})^n \; \& \; h \in \hbox{cyl}(B)\\
&\quad\ \ \,\&\ A \cap
\hbox{cyl}(B) = \phi)\\[.5pc]
& \leftrightarrow (\exists n \geq 1)\; (\exists B)\;(B \;\mbox{
is a }\; \Delta^1_1 \; \hbox{subset of} \; (\omega^{\omega})^n\\
&\quad\ \ \,\& \; h \;\in \;\hbox{cyl}(B)\; \& \; A \;\cap
\;\hbox{cyl}(B) = \phi).
\end{align*}
To prove the previous implication $\rightarrow$, let $B$ be a
$\Sigma^1_1$ subset of $(\omega^{\omega})^n$ such that $h \in
\hbox{cyl}(B)$ and $A \cap \hbox{cyl}(B) = \phi$. But then $p_n(A)
\cap B = \phi$. Since $p_n(A)$ is $\Sigma^1_1$, it follows from
Kleene's separation theorem that there is a $\Delta^1_1$ subset
$B'$ of $(\omega^{\omega})^n$ such that $B \subseteq B'$ and $B'
\cap p_n(A) = \phi$. Hence $h \in \hbox{cyl}(B')$ and $A \cap
\hbox{cyl}(B') = \phi$, which establishes $\rightarrow$.
Consequently,
\begin{equation*}
h \notin \hbox{cl}_1(A) \leftrightarrow (\exists \; m)(m \in W_0 \
\& \ C(m,h) \ \& \ C_m \cap A  = \phi).
\end{equation*}
So $\neg \hbox{cl}_1(A)$ is $\Pi^1_1$.\hfill $\Box$
\end{proof}

\begin{lem}
Assume {\rm ({\bf O})}. If $A$ is a $\Pi_1$ subset of $H${\rm ,}
then $A \Delta {\rm sur}_1(A)$ is $T^\infty$-meager.
\end{lem}

\begin{proof}
Choose subsets $B_n$ of $(\omega^{\omega})^n, \; n \geq 1$, such
that
\begin{equation*}
A = H - \cup_{n \geq 1} \hbox{cyl}(B_n).
\end{equation*}
Then
\begin{align*}
\hbox{sur}_1(A) - A &= \hbox{sur}_1(A) \cap \cup_{n \geq 1}
\hbox{cyl}(B_n)\nonumber\\[.5pc]
&\subseteq \cup_{n \geq 1} ([\hbox{sur}_1(A) \cap
\hbox{cyl}(\hbox{sur}^n(B_n))]\\[.2pc]
&\quad\,\cup[\hbox{cyl}(B_n) - \hbox{cyl}(\hbox{sur}^n(B_n))]).
\end{align*}
Now
\begin{equation*}
\hbox{cyl}(B_n) - \hbox{cyl}(\hbox{sur}^n(B_n))  = \hbox{cyl}(B_n
- \hbox{sur}^n(B_n)).
\end{equation*}
The set on the right of  the above equality is $T^\infty$-meager
by virtue of Lemma~2.13 in [L2].  We will now prove that
$\hbox{sur}_1(A) \cap \hbox{cyl}(\hbox{sur}^n(B_n))$ is
$T^\infty$-nowhere dense. Note that $\hbox{sur}_1(A) \cap
\hbox{cyl}(\hbox{sur}^n(B_n))$ is $T_1$-closed, hence
$T^\infty$-closed. Now let $A'$ be a $\Sigma^1_1$ set contained in
$\hbox{sur}_1(A) \cap \hbox{cyl}(\hbox{sur}^n(B_n))$. Then
\begin{equation*}
\hbox{cyl}(p_n(A')) \subseteq \hbox{cyl}(\hbox{sur}^n(B_n)).
\end{equation*}
Hence
\begin{align*}
A \cap \hbox{cyl}(p_n(A')) & \subseteq \hbox{cyl}(\hbox{sur}^n(B_n)) -
\hbox{cyl}(B_n)\\[.5pc]
& = \hbox{cyl}(\hbox{sur}^n(B_n) - B_n).
\end{align*}
Consequently, by virtue of Lemma 2.2 and Lemma 2.13 in [L2], $A
\cap \hbox{cyl}(p_n(A'))$ is $T^\infty$-meager. Since
$\hbox{cyl}(p_n(A'))$ is $T_1$-open, it follows that
$\hbox{cyl}(p_n(A')) \subseteq \hbox{cosur}_1(A)$. Hence $A'$ is
empty because $A'$ is also contained in $\hbox{sur}_1(A)$. Thus
$\hbox{sur}_1(A) \cap \hbox{cyl}(\hbox{sur}^n(B_n))$ is
$T^\infty$-nowhere dense. It follow from (1) that
$\hbox{sur}_1(A)-A$ is $T^\infty$-meager. Since
$A-\hbox{sur}_1(A)$ is easily seen to be $T^\infty$-meager, we are
done. \hfill $\Box$
\end{proof}

\begin{lem}
If $A$ and $B$ are $\Sigma^1_1$ subsets of $H$ such that $A$ can
be separated from $B$ by a $\Sigma_1$ set{\rm ,} then $A \cap {\rm
cl}_1(B) = \phi.$
\end{lem}

\begin{proof}
Suppose $D$ is a $\Pi_1$ subset of $H$ such that $A \cap D = \phi$
and $B \subseteq D$. Hence, by Lemma 3.2, $B-\hbox{sur}_1(D)$ is
$T^\infty$-meager. But $B - \hbox{sur}_1(D)$ is $T^\infty$-open,
so $B \subseteq \hbox{sur}_1(D)$.

Since $\hbox{sur}_1(D)$ is $T_1$-closed, $\hbox{cl}_1(B) \subseteq
\hbox{sur}_1(D)$. Now $A \cap \hbox{sur}_1(D)$ is
$T^\infty$-meager, so $A \cap \hbox{cl}_1(B)$ is
$T^\infty$-meager. By Lemma 3.1, $A \cap \hbox{cl}_1(B)$ is
$\Sigma^1_1$, hence $A \cap \hbox{cl}_1(B)$ must be empty.\hfill
$\Box$
\end{proof}

\begin{lem}
If $A$ and $B$ are $\Sigma^1_1$ subsets of $H$ such that $A \cap
{\rm cl}_1(B) = \phi${\rm ,} then $A$ can be separated from $B$ by
a $\Sigma^*_1$ set.
\end{lem}

\begin{proof}
Define
\begin{equation*}
P(h,n) \leftrightarrow h \notin A \vee (n \in W_0 \ \& \ C(n,h) \
\& \ C_n \cap B = \phi).
\end{equation*}
Then $P$ is $\Pi^1_1$ and $(\forall h)(\exists n)P(h,n)$. By
Kreisel's selection theorem $[\hbox{Mo}]$, there is a
$\Delta^1_1$-recursive function $f\hbox{\rm :}\ H \rightarrow
\omega$ such that $(\forall h)P(h,f(h))$. Let
\begin{equation*}
D = \{ n \in \omega\hbox{\rm :}\ n \in W_0\ \&\ C_n \cap B = \phi
\}.
\end{equation*}
Then $D$ is $\Pi^1_1$ and $f(A) \subseteq D$. Since $f(A)$ is
$\Sigma^1_1$, there is a $\Delta^1_1$ set $E \subseteq \omega$
such that $f(A) \subseteq E \subseteq D$. Let
\begin{equation*}
R(h,n) \leftrightarrow n \in E \ \& \ C(n,h),
\end{equation*}
Then $R$ is $\Delta^1_1$, because if
\begin{equation*}
R'(h,n) \leftrightarrow n \in E \ \& \ \neg C(n,h),
\end{equation*}
then both $R$ and $R'$ are $\Pi^1_1$, $R \cap R' = \phi$ and $R
\cup R' = H \times E$. Set
\begin{equation*}
G_n = \{ h\hbox{\rm :}\ R(h,n)\}, \; n \in \omega.
\end{equation*}
Then $\cup_{n \geq 0} G_n$ is a $\Sigma^*_1$ set which separates
$A$ from $B$.\hfill $\Box$
\end{proof}

Lemmas 3.2, 3.3 and 3.4 establish Theorem 1.1 for $\xi = 1$.

\section{Proof of Theorem 1.1}

The proof of Theorem 1.1 is by induction on $\xi$. So we fix $\xi
> 1$ and assume Theorem 1.1 is true for all $\eta < \xi$. Lemmas
3.1--3.4 can be formulated and proved at level $\xi$, thereby
completing the proof of Theorem 1.1 at level $\xi$. We omit the
proofs because they are exactly like the proofs of Lemmas 7, 8, 9
and Theorem B in [L1].

We observe that the inductive hypothesis that Theorem 1.1 hold at
all levels $\eta < \xi$ is by itself not sufficiently strong to
prove the analogue of Lemma 3.2 at level $\xi$ and hence the
theorem itself at that level. For this we need that analogues of
Lemma 3.2 hold at all levels $\eta < \xi$. It is at this point in
the proof that assumption ({\bf O}) is needed to ensure that Lemma
3.2 hold at level $\xi=1$,  the higher levels of Lemma 3.2 then
being proved by inducting up from the base level.

\section{Concluding remarks}
\setcounter{defin}{0}

For $\alpha  \in  \omega^{\omega}$, we now consider the following
statement of set theory:\vspace{.5pc}

\noindent $(\alpha)$ Every subset of $\omega^{\omega}$ has the
Baire property with respect to the topology whose base is the
pointclass of $\Sigma^1_1(\alpha)$ sets.\vspace{.5pc}

It is straightforward to relativize Theorem 1.1 to $\alpha$ under
the assumption that $(\alpha)$ holds. The next result is provable
in $\hbox{ZF}+\hbox{DC}+(\forall\alpha)((\alpha))$.

\begin{theo}[\!]
Let $X$ be an uncountable Polish space and let $H = X^N$. Then{\rm
,} for $1 \leq \xi < \omega_1${\rm ,}
\begin{equation*}
\overline{\Sigma}_\xi = \Sigma_\xi \cap {\cal B}.
\end{equation*}
\end{theo}

Under the assumption that there is an inaccessible cardinal,
Solovay [S] proved that $\hbox{ZF}+\hbox{DC}$ holds in the
L\'{e}vy--Solovay model. Furthermore, it was observed  by Louveau
(p.43 of [L2]) that the statement $(\forall \alpha)((\alpha))$
holds as well in the  model.

Whether Theorem 5.1 is provable in ZFC remains an open problem.
Indeed, we do not have an answer to the problem even when $\xi=2$.

It is not difficult to prove that the axiom of determinacy implies
$(\forall\alpha)((\alpha))$ so that Theorem 5.1 is provable in
$\hbox{ZF}+\hbox{AD}$ (see [Mo]). On the other hand, the axiom of
choice implies $\neg({\bf O})$ in ZF.

\section*{Appendix}

We will now prove Theorem 1.1 for $\xi=1$  without assuming ({\bf
O}). In view of Lemma 3.4, it will suffice to prove that $A\cap
\hbox{cl}_1(B)=\phi$. Define
\begin{equation*}
P(h,n) \leftrightarrow (n \geq 1)\ \& \ (\exists h')(p_n(h)h' \in
B),
\end{equation*}
where $p_n(h)h'$ is the catenation of $p_n(h)$ and $h'$. Note that
$P$ is $\Sigma^1_1$. Let
\begin{equation*}
h \in \bar{B} \leftrightarrow (\forall n \geq 1)P(h,n),
\end{equation*}
so that $\bar{B}$ is the closure of $B$ in the product of discrete
topologies on $H$. Consequently, $\bar{B} \subseteq H - A.$ Define
\begin{equation*}
Q(h,n) \leftrightarrow (n \geq 1) \ \& \ (\neg P(h,n) \vee h
\notin A).
\end{equation*}
Then $Q$ is clearly $\Pi^1_1$ and $(\exists n)Q(h,n)$. So there is
a $\Delta^1_1$-recursive function $f\hbox{\rm :}\ H \rightarrow
\omega$ such that $(\forall h) Q(h,f(h)).$ Let\pagebreak
\begin{equation*}
S(h,n) \leftrightarrow (n \geq 1) \ \& \ (f(h) \neq n \; \vee h
\notin A).
\end{equation*}
{\it Claim.}

\begin{enumerate}
\renewcommand\labelenumi{(\roman{enumi})}
\leftskip .3pc
\item $S \ \mbox{is} \ \Pi^1_1$,

\item $(\forall h)(\forall n \geq 1) (P(h,n) \rightarrow
S(h,n))$,

\item $h \notin A \leftrightarrow (\forall n \geq 1) S(h,n)$.
\end{enumerate}

To see (ii), assume $P(h,n)$. Then we must have $h \in A
\rightarrow f(h) \neq n$. Hence $S(h,n)$. For (iii), suppose $h
\notin A$. Clearly, then $(\forall  n \geq 1)S(h,n)$. Suppose now
that $h \in A$. Then there is $n$ such that $f(h) = n$, hence
$\neg S(h,n)$. (iii) now follows.

Now turn each $S_n$ into a cylinder set as follows. Define
\begin{equation*}
R(h,n) \leftrightarrow (\forall h')S(p_n(h)h', n),
\end{equation*}
so $R$ is $\Pi^1_1$. Note that $P_n$ and $R_n$ are cylinder sets,
that is,
\begin{equation*}
P(h,n) \ \& \ p_n(h) = p_n(h') \rightarrow P(h',n)
\end{equation*}
and
\begin{equation*}
R(h,n)\ \&\ p_n(h) = p_n(h') \rightarrow R(h',n).
\end{equation*}
{\it Claim.} $(\forall h) (\forall n) (P(h,n) \rightarrow
R(h,n)).$\vspace{.5pc}

So suppose $P(h,n)$. Then, for every $h'$, $P(p_n(h)h',n)$, hence
$S(p_n(h)h',n)$, so $R(h,n)$.

To complete the proof,  let $h \in A$. Then there is $n \geq 1$
such that $\neg S(h,n)$, hence $\neg R(h,n)$. Now $\neg R_n$ is
$\Sigma^1_1$ and $\Pi_0$ because $R_n$ is a cylinder set.
Moreover, $\neg R_n \cap B = \phi$ because $\neg R_n \subseteq
\neg P_n$ and $\neg P_n \cap B = \phi$. Hence $\neg R_n$ is a
$T_1$-open set containing $h$ and disjoint from $B$. So $h \notin
\hbox{cl}_1(B)$.

\section*{Acknowledgement}

The authors would like to thank the referee for making a number of
helpful suggestions.


\begin{thebibliography}{999}
\bibitem[Ba]{Ba} Barua~R, On Borel hierarchies of countable products of Polish
spaces, {\it Real Analysis Exchange} {\bf 16} (1990--91) 60--66

\bibitem[L1]{L1} Louveau~A, A separation theorem for $\Sigma^1_1$ sets, {\it Trans.
Am. Math. Soc.} {\bf 260} (1980) 363--378

\bibitem[L2]{L2} Louveau~A, Ensembles analytiques et bor\'{e}liens dans les
espaces produits, {\it Ast\'{e}risque} {\bf 78} (1980) 1--87

\bibitem[Ma]{Ma} Maitra~A, Pestien~V and Ramakrishnan~S, Domination by Borel
stopping times and some separation properties, {\it Fund. Math.}
{\bf 135} (1990) 189--201

\bibitem[Mo]{Mo} Moschovakis~Y~N, Descriptive Set Theory
(Amsterdam: North-Holland) (1980)

\bibitem[S]{S} Solovay~R~M, A model of set theory in which every set of
reals is Lebesgue measurable, {\it Ann. Math.} {\bf 92} (1970)
1--56
\end{thebibliography}
\end{document}